\documentclass[11pt]{article}
\usepackage{graphicx}
\usepackage[dvips]{epsfig}
\def\institute#1{\large #1}

\begin{document}
\title{Ecoepidemics with a nonlinear disease incidence}
\author{Ezio Venturino\\[2mm] \institute{Dipartimento di Matematica ``Giuseppe Peano'',\\
Universit\`{a} di Torino, Italy.}}
\date{}
\maketitle
\thispagestyle{plain}
\begin{abstract}
We present two new models for interacting populations subject to a transmissible
disease. The novelty lies in the assumption that herd behavior influences the disease
incidence, rather than the demographic description of the interactions, as in previous
related similar models. As it is already known from other ecoepidemiological situations,
the epidemics may affect the system demographic outcomes.
\end{abstract}
\section{Introduction}

Ecoepidemiology studies the influence of diseases among interacting populations.
This rather new field of research started about a quarter of a century ago, with
investigations merging diseases in demographic models in different contexts,
\cite{HF,BC,V94}. For a brief overview of the progress up to a few years ago,
see Chapter 7 of \cite{MPV}.

Much more recently, a novel idea for modeling herd behavior has been introduced, \cite{AV,APV},
and further explored in \cite{CV}. It is essentially based on the observation that individuals
gathering in groups can generally be attacked by predators on the outskirts of the
territory that they occupy, that is proportional to their size. The population occupies
thus a two-dimensional manifold, while its boundary represents a one-dimensional one.
The former is directly related to the population size, while the latter instead must
be then related to its square root.

Rather then pursuing 
this idea in various circumstances in ecological situations, \cite{APV},
following the idea of extending these demographic remarks to ecoepidemic situations, \cite{EV11,CV},
we want to exploit it here still in the ecoepidemic realm, but considering gatherings
of infected individuals.
In fact, we use the basic ideas on herd behavior, \cite{APV},
in a different context. Specifically, in \cite{EV11} as well as in \cite{CV,CV13}, 
it is still the demographic part of the model that is modeled according to herd behavior.
Here, however, we assume that the infected lump together. Their herd size grows due to
the arrival of new individuals
only through the susceptible contacts, wandering about independently of each other,
with the infectious individuals positioned
on the outskirts of the
infected bunch. This may be plausible in the context of predators and prey interactions,
because infected individuals, in general weaker
and slower, may gather together. Hence, their new possible recruits would arrive precisely
through the above mechanism. 
As stated above, assuming thus
that these populations occupy a certain portion of ground, predation occurs on the border,
i.e. the perimeter, of the
lump of infected individuals, therefore it is expressed via a square root term of their size.
 
The paper is organised as follows. We briefly summarize the results on the classical
predator-prey reference model in the next Section, then provide some basic information
on the ecoepidemic models we want to introduce. In Section \ref{sec:1} we consider the
case of infected prey that are harmless for predators. Section \ref{sec:2} contains
the particular case in which the predators recognize and avoid infected prey.
In Section \ref{sec:3} we present the model for which infected are toxic for the predators.
Results on the boundedness of the systems' trajectories are derived in the next Section
and a final summary of the results concludes the paper.

\section{The classical predator-prey reference model}

Let us consider the Lotka-Volterra predator-prey model with logistic correction for the prey $Q$.
The predators $P$ are assumed to be specialists,
so that in the absence of $Q$ they would starve to death.
The model reads
\begin{eqnarray}\label{Mod_0}
\frac{dP}{dt} &=& -mP+aPQ, \\ \nonumber 
\frac{dQ}{dt} &=& rQ\left( 1-\frac{Q}{K}\right)-aPQ.
\end{eqnarray}
The parameters are defined as follows:
$m$ represents the predator's mortality rate,
$a$ is the predator's hunting rate,
$r$ is the prey reproduction rate,
$K$ is the prey carrying capacity.

This system dynamics is well known. There are only two possibly stable equilibria, since the origin
is always unstable. Between the predator-free point $E_1^c = (0,K)$
and coexistence $\widehat E_*^c = (P_*^c,Q_*^c)$
there is a transcritical bifurcation.
Whenever the predators' mortality rate falls below the threshold
\begin{equation}\label{m*}
m^*=aK
\end{equation} 
the predators invade the environment permanently, as $E_*^c$ becomes feasible and
is unconditionally stable, while instead $E_1^c$ loses its stability.
In each case, the only possible equilibrium is globally asymptotically stable.

\section{Background on the ecoepidemic situation}

In the ecoepidemic approach the main difference with the previous demographic model consists
in the fact that
the prey population is divided among infected $I$, that gather in herds, and susceptibles $S$.
We also assume the presence in the ecosystem of a third population
$P$ that can predate on $S$ as well as possibly on the $I$ lumped together. With this we mean
that predation occurs always on the border of the lump of $I$'s, as in other herd
behavior systems, \cite{APV}. But the effect
of predation on the infected can have different outcomes for the predators.

Three cases will be considered, as far as the behavior of $P$ with respect to $I$ is
concerned: they do not recognize the infected, but their predation leaves the $P$'s unaffected,
they recognize the $I$'s and avoid them, or finally the
$P$'s predate the $I$'s and the latter harm the predators.
We stress here once again that the novelty of this model is in using the herd behavior in the
epidemiological terms, especially in contrast with what was done in \cite{CV,CV13},
where the major issue was on the infected behavior, but still considered from the demographic
behavior point of view.

The group gathering behavior is modeled as indicated above, via the square root of
the infected population density $\sqrt I$. Disease transmission occurs through
contacts among the infected lying on the boundary of the herd with the susceptibles.
Assuming homogeneous mixing among these classes of individuals, the corresponding
(modified) ``mass action'' term assumes the form $S\sqrt I$.
Thus,
our nonlinear disease incidence model could be regarded as a
particular case of the $S^{\alpha}I^{\gamma}$ incidence, which has been proposed
among other epidemiological population interaction possibilities, see \cite{Liu1986}.
f
The following are further general assumptions for all the three models considered here:
the lump of $I$'s does not reproduce, it can grow only by recruiting newcomers from the class $S$.
They are also too weak to exert any intraspecific competition on the healthy individuals $S$,
nor feel their pressure for the search for resources, since they do not reproduce.
We also assume that the predators
do not have other food sources, being specialists.
With respect to the loose population $S$ the encounters with $P$ are on a one to
one basis, i.e. they are expressed via the usual mass action law.

\section{Infected are harmless for predators.}\label{sec:1} 

The first model we investigate is the following one
\begin{eqnarray}\label{Mod_2}
\frac{dP}{dt} &=& -mP+aPS+bP\sqrt{I},\\ \nonumber 
\frac{dS}{dt} &=& -\beta S\sqrt{I}+rS\left( 1-\frac{S}{K}\right)-aPS, \\ \nonumber 
\frac{dI}{dt} &=& -\mu I+\beta S\sqrt{I}-bP\sqrt{I} .
\end{eqnarray}
We define the meaning of all the parameters, because although we use the same notation
as in (\ref{Mod_0}), the interpretation of the parameters common to both (\ref{Mod_0})
and (\ref{Mod_2}) is at times slightly different.
Here 
$m$ represents the predator's mortality rate,
$a$ is the predator's hunting rate on healthy prey,
$b$ is the predation rate on the infected herds of prey,
$\beta$ is the disease incidence rate,
$r$ is the healthy prey reproduction rate,
$K$ is the carrying capacity of healthy prey,
$\mu$ is the natural plus disease-related mortality rate of infected individuals.

In view of the fact that the prey modeled by the $I$'s lump together,
predation on them is 
exerted only on the outer boundary of their herd, which is expressed by
the square root term in the above first and third equations (\ref{Mod_2}).
We need to redefine the dependent variables, to avoid a possible singularity in the Jacobian
when $I$ vanishes.
Singularity removal can be performed by defining $U=\sqrt{I}$. It leads to
\begin{eqnarray}\label{caso2bis}
\frac{dP}{dt} &=& P\left( -m+aS+bU\right) ,\\ \nonumber
\frac{dS}{dt} &=& S\left[ -\beta U+r\left( 1-\frac{S}{K}\right)-aP\right] ,\\ \nonumber
\frac{dU}{dt} &=& \frac{1}{2} \left(-\mu U+\beta S-bP\right).
\end{eqnarray}
The Jacobian of (\ref{caso2bis}) is
\begin{eqnarray}\label{Jac_caso2bis}
J^h=\left(
\begin{array}{ccc}
-m +a S +bU & aP & bP\\
-aS & -\beta U +r-2\frac rK S -aP & -\beta S\\
-\frac 12 b & \frac 12 \beta & -\frac 12 \mu
\end{array}
\right).
\end{eqnarray}

\subsection{Equilibria and their analysis}

The possible equilibria are the points
$E_0$, namely the system disappearance, the predator-free point
$$
E_1= \left(0, \ rK\frac{\mu}{\beta ^2 K+r\mu}, \ rK\frac{\beta}{\beta ^2 K+r\mu}\right) ,
$$
and the coexistence
equilibrium $\widehat E_*$ with population values
\begin{eqnarray*}
\widehat P_*=\frac{bKr\beta +aKr\mu -mr\mu-Km\beta ^2}{a^2K\mu +b^2r}, \quad
\widehat S_*=K\frac{am\mu +b^2 r-bm\beta }{a^2K\mu +b^2r}, \\
\widehat U_*=\frac{bmr+aKm\beta -abKr}{a^2K\mu +b^2r}.
\end{eqnarray*}
Feasibility conditions for $\widehat E_*$ are 
\begin{equation}\label{E*2_feas}
Kr(b\beta +a\mu) \ge m(r\mu+K\beta ^2), \quad
am\mu +b^2 r\ge bm\beta , \quad
m(br+aK\beta )>abKr.
\end{equation}

The origin is unstable, since the eigenvalues of $J^h$ are $-m$, $r$, $-\frac 12 \mu$.

At $E_1$ the stability condition is regulated by the very first eigenvalue, 
\begin{equation}\label{E1_2_stab}
aS_1 + b U_1 \equiv m^{\ddagger} <m,
\end{equation}
since the remaining ones come from a 2 by 2 submatrix $J^h_2$ for which the Routh-Hurwitz conditions
hold unconditionally, since they become
\begin{equation}\label{RH_E1}
-tr(J^h_2)=\frac rK S_1 + \frac 12 \mu>0, \quad
\det(J^h_2)=\frac rK S_1 +\frac 12 \beta^2 S_1 > 0.
\end{equation}

Note that (\ref{E1_2_stab}) is the opposite condition of the first inequality for
the feasibility of $\widehat E_*$, (\ref{E*2_feas}), so
that when the other two feasibility conditions (\ref{E*2_feas}) hold,
we have a transcritical bifurcation for which $\widehat E_*$ emanates from $E_1$.
It is clearly seen also that no Hopf bifurcation can arise here, in view of the strict inequality for
the trace.

The coexistence equilibrium $\widehat E_*$ is always stable, whenever feasible, since the characteristic
equation (\ref{char}) is the cubic
\begin{equation}\label{char}
\sum _{k=0}^3 a_k \lambda ^k=0,
\end{equation}
with the coefficients
\begin{eqnarray}\label{coeffs}
a_0=\frac 12 \left( a^2 \mu + b^2 \frac rK \right) S_* P_* > 0, \quad 
a_2=\frac rK S_* + \frac 12 \mu  > 0,\\ \nonumber
a_1=a^2 S_* P_* + \frac 12 \left[ b^2 P_3 + \left( \frac rK \mu + \beta^2 \right) S_* \right] > 0.
\end{eqnarray}
In fact, also the last Routh-Hurwitz conditions holds unconditionally
\begin{equation}\label{RH_E2}
a_2a_1-a_0=\frac rK S_*^2 \left( a^2 P_* + \frac {r\mu} {2K} + \frac 12 \beta^2 \right)
+ \frac {\mu}4 \left( \frac {r\mu}K S_* + \beta^2 S_* + b^2  P_*\right) >0 
\end{equation}
and strictly, thus preventing also possible Hopf bifurcations.

\section{No predation on infected prey}\label{sec:2} 

We briefly examine here the particular case of (\ref{Mod_2}) in which $b=0$, i.e. the infected
are recognized and completely
disregarded by the predators. The system with no singularity
and its Jacobian are obtained just as particular cases of (\ref{caso2bis})
and (\ref{Jac_caso2bis}), setting in them $b=0$.

The possible equilibria are again all the points found earlier, namely the origin and the predator-free
equilibrium $E_1$ with the very same population values. Both these equilibria
are both always feasible. We also find coexistence, which now simplifies to
$$
E^*= \left( \frac{arK\mu-m\beta ^2 K-mr\mu}{a^2K\mu}, \ \frac{m}{a}, \ \frac{\beta m}{\mu a}\right) .
$$
It is feasible only if
\begin{equation}\label{E1_feas}
m< m^{\dagger} \equiv \frac {aKr\mu}{K\beta^2+r\mu}.
\end{equation} 
This condition specifies that the predator's mortality must fall below a certain critical threshold.
Note that $m^{\dagger}$ coincides with $m^{\ddagger}$ when the latter is evaluated for $b=0$.

The origin $E_0$ retains its unconditional instability, in view of the very same eigenvalues we found for
(\ref{caso2bis}), namely $-m$, $r$, $-\frac 12 \mu$.

One eigenvalue of $E_1$ is now $-m+aS_1$ giving the stability condition,
as the last two conditions (\ref{E*2_feas}) now are trivially satisfied:
\begin{equation}\label{E1_stab}
m^{\dagger} <m,
\end{equation}
again a particular case of what we found for (\ref{caso2bis}). For the remaining ones again
the Routh-Hurwitz conditions hold unconditionally. Indeed
the remaining ones come from the 2 by 2 submatrix $\widetilde J_2^h=J_2^h$ for which the Routh-Hurwitz conditions
hold unconditionally (\ref{RH_E1}). Again, no Hopf bifurcation can arise here as well
and for $m=m^{\dagger}$ there is a transcritical bifurcation for which $E^*$ arises from $E_1$.

At $E^*$ the characteristic equation is the cubic
(\ref{char})
with coefficients that are obtained from (\ref{coeffs}) by setting $b=0$. Therefore
since all these coefficients are strictly positive and also the third
Routh-Hurwitz stability condition holds, whenever feasible, the coexistence
equilibrium is unconditionally stable. Also, in view of the above strict inequality
in (\ref{RH_E2}), no Hopf bifurcations can arise here as well.

\section{The case of toxic infected.}\label{sec:3} 

In this case we assume that the infected prey are harmful for the predators when
they come in contact. The model becomes then
\begin{eqnarray}\label{Mod_3}
\frac{dP}{dt} &=& -mP+aPS-bP\sqrt{I},\\ \nonumber
\frac{dS}{dt} &=& -\beta S\sqrt{I}+rS\left( 1-\frac{S}{K}\right)-aPS ,\\ \nonumber
\frac{dI}{dt} &=& -\mu I+\beta S\sqrt{I}-bP\sqrt{I}.
\end{eqnarray}
Again, all the parameter retain their meaning from (\ref{caso2bis}), but note the change
in sign in the last term of the first equation.
Once again, the system with the removed singularity becomes
\begin{eqnarray}\label{caso3bis}
\frac{dP}{dt} &=& P\left( -m+aS-bU\right) ,\\ \nonumber
\frac{dS}{dt} &=& S\left[ -\beta U+r\left( 1-\frac{S}{K}\right)-aP\right] ,\\ \nonumber
\frac{dU}{dt} &=& \frac{1}{2} \left(-\mu U+\beta S-bP\right).
\end{eqnarray}

The Jacobian of (\ref{caso3bis}) is
$$
J^t=\left(
\begin{array}{ccc}
-m +a S -bU & aP & -bP\\
-aS & -\beta U +r-2\frac rK S -aP & -\beta S\\
-\frac 12 b & \frac 12 \beta & -\frac 12 \mu
\end{array}
\right).
$$

\subsection{Equilibria and their analysis}

Again the origin $E_0$ and the predator-free equilibria $E_1$
are unaltered from the previous case
(\ref{caso2bis})
and are therefore always feasible. Coexistence $\widetilde E_*$
settles instead at the following population values
\begin{eqnarray*}
\widetilde P_*=\frac{aKr \mu - mr\mu -Km\beta ^2 -bKr\beta}{a^2K\mu - 2abK\beta -b^2 r }, \quad
\widetilde S_*=\frac{aKm\mu -bKm\beta -b^2Kr }{a^2K\mu -2abK\beta -b^2 r } \\
\widetilde U_*=\frac{aKm\beta +bmr -abKr }{a^2K\mu -2abK\beta -b^2 r }.
\end{eqnarray*}
Feasibility conditions are either one of these sets of inequalities
\begin{eqnarray}\label{E*tilde_f1}
2abK\beta +b^2 r \ge a^2K\mu, \quad
mr\mu+Km\beta ^2+bKr\beta \ge aKr \mu, \\ \nonumber
bKm\beta + b^2Kr \ge aKm\mu, \quad
abKr \ge bmr+aKm\beta ;
\end{eqnarray}
or
\begin{eqnarray}\label{E*tilde_f2}
2abK\beta +b^2 r \le a^2K\mu, \quad
mr\mu+Km\beta ^2+bKr\beta \le aKr \mu, \\ \nonumber
bKm\beta + b^2Kr \le aKm\mu, \quad
abKr \le bmr+aKm\beta .
\end{eqnarray}

Stability of $E_0$ is once again unchanged, the eigenvalues are still the same, $-m$, $r$, $-\frac 12 \mu$.
For $E_1$ we find again the very same condition (\ref{E1_2_stab}), as the remaining analysis on the 2 by 2
submatrix carries out unaltered.

\begin{figure}[ht]
\centering
\includegraphics[width=8.5cm]{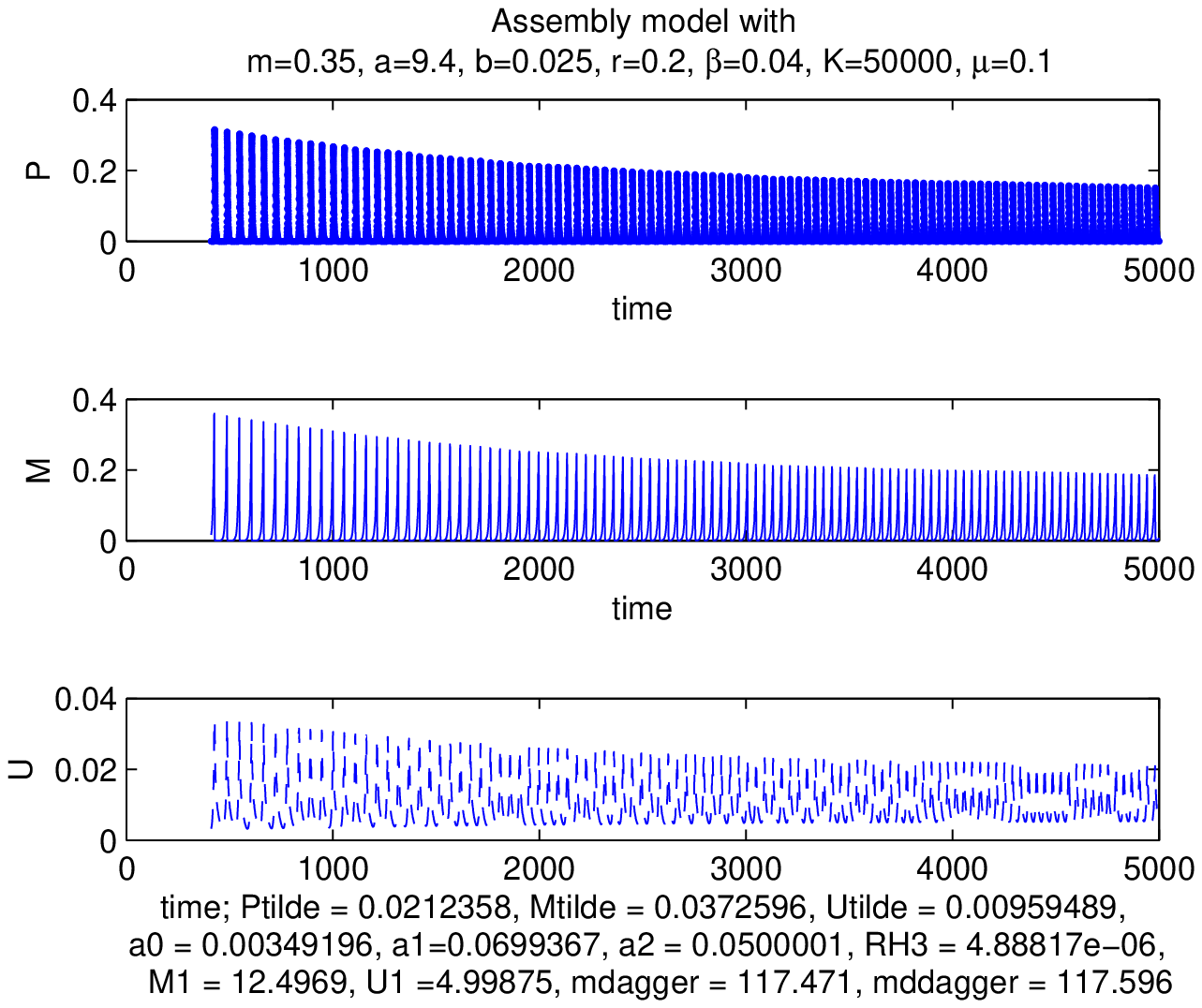}
\caption{For the system (\ref{caso3bis}), decaying oscillations involving all three subpopulations
arise for the parameter values
$m=0.35$, $a=9.4$, $b=0.025$, $r=0.2$, $\beta=0.04$, $K=50000$, $\mu=0.1$, dampened toward the equilibrium
$\widetilde E_* =(0.0212, 0.0373, 0.0096)$.
The Routh-Hurwitz conditions do however hold:
$a_0=0.0035$, $a_1=0.0699$, $a_2=0.0500$ and $a_2a_1-a_0=8.607 \times 10^{-6}$.
Top to bottom the populations $P$, $S$, $U$, as functions of time.}
\label{fig:d_oscill_PMU}
\end{figure}

At coexistence instead relevant changes occur, as the cubic (\ref{char}) has here the coefficients
\begin{eqnarray*}
a_0=\frac 12 \left( a^2 \mu - b^2 \frac rK -2ab\beta \right) \widetilde S_* \widetilde P_*, \quad 
a_2=\frac rK \widetilde S_* + \frac 12 \mu , \\
a_1=a^2 \widetilde S_* \widetilde P_* + \frac 12 \left[ \left( \frac rK \mu + \beta^2 \right) \widetilde S_*
- b^2 \widetilde P_*\right].
\end{eqnarray*}
They now must all be imposed to be positive. Note that $a_0>0$ is incompatible with the first feasibility
condition for $\widetilde E_*$ (\ref{E*tilde_f1}).
Also the condition $a_2a_1-a_0>0$ must be imposed, which now becomes
\begin{eqnarray}\label{E*tilde_HB}
\frac rK \widetilde S_*^2 \left( a^2 \widetilde P_* +
\frac r{2K} \mu + \frac 12 \beta^2 \right)
+ \frac 14 \mu \widetilde S_* \left( \frac rK \mu + \beta^2 \right)
+ab \beta \widetilde S_* \widetilde P_* 
> \frac 14 \mu b^2  \widetilde P_*.
\end{eqnarray}
It would therefore in principle be possible that Hopf bifurcations in this case could arise.
However extended simulations attempting to violate this conditions were not successful. We have
been able only to make it almost an equality, but never to reverse the above inequality
(\ref{E*tilde_HB}), see Fig. \ref{fig:d_oscill_PMU}.

\begin{figure}[ht]
\centering
\includegraphics[width=8.5cm]{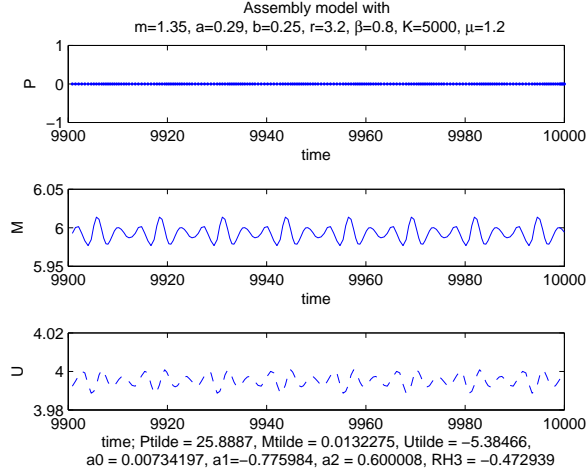}
\caption{For the system (\ref{caso3bis}),
the prey subpopulations can thrive together, in absence of predators, also via tiny persistent
oscillations, here obtained with the parameter values
$m=1.35$, $a=0.29$, $b=0.25$, $r=3.2$, $\beta=0.8$, $K=5000$, $\mu=1.2$.
The coexistence equilibrium in this case is unfeasible, $\widetilde E_* =(25.889, 0.013, -5.385)$.
In this case $E_1 =(0, 5.993, 3.995)$ and $m<m^{\ddagger}=2.737$, showing its instability, compare (\ref{E1_2_stab}).
Note that the oscillations shown are indeed around this predator-free equilibrium point.
Top to bottom the populations $P$, $S$, $U$, as functions of time.}
\label{fig:oscill_MU}
\end{figure}

We conjecture therefore that also in this case the coexistence equilibrium does not lead to
Hopf bifurcations. Instead, mainly by rendering $a_0$ negative, we can destabilize the coexistence equilibrium.
This in turn takes the system trajectories either to the predator-free equilibrium $E_1$,
when it is stable, namely for (\ref{E1_2_stab}), or to limit cycles around it, see Fig. \ref{fig:oscill_MU}.

\begin{figure}[ht]
\centering
\includegraphics[width=8.5cm]{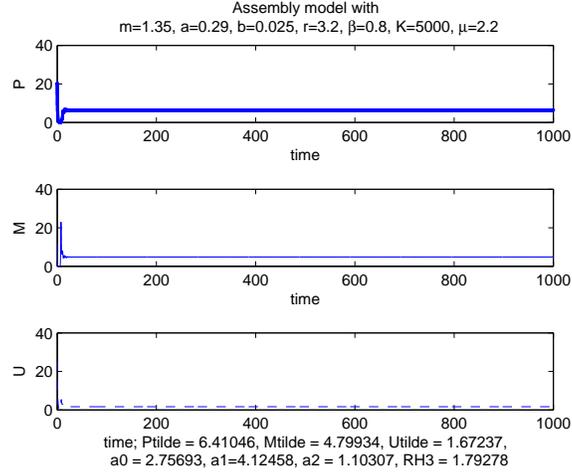}
\caption{For the system (\ref{caso3bis}),
the coexistence equilibrium can be stably achieved for the parameter values
$m=1.35$, $a=0.29$, $b=0.025$, $r=3.2$, $\beta=0.8$, $K=5000$, $\mu=2.2$, at the level
$\widetilde E_* =(6.410, 4.799, 1.672)$. The Routh-Hurwitz conditions hold, since
$a_0=2.667$, $a_1=4.125$, $a_2=1.103$ and $a_2a_1-a_0=1.882$. 
Top to bottom the populations $P$, $S$, $U$, as functions of time.}
\label{fig:coex}
\end{figure}

\section{Boundedness}\label{sec:4} 

The finiteness of the trajectories can be shown for all three original models together as follows.
Let $T=P+S+I$, by adding the differential equations
it is then easy to show that for (\ref{Mod_2}) and (\ref{Mod_3}) the
following inequality holds:
$$
T\le -mP + rS -\frac rK S^2 -\mu I.
$$
Taking now an arbitrary $0<q < \min \{\mu, m\}= M$, we find
$$
\frac {dT}{dt} + q T \le (r+q) S - \frac rK S^2 + (q-M) (P+I) \le \Psi,
$$
since $q-M<0$ and
where $\Psi$ denotes the height of the vertex of the parabola in $S$ on the left hand side, for which
$$
\Psi = \frac K{4r} (r+q)^2.
$$
It follows then that the solutions of the above differential inequality must lie below those of
$$
\frac {dT}{dt}=\Psi - qT,
$$
i.e.
$$
T(t) = \frac {\Psi}q \left( 1 - e^{-qt} \right) + T(0) e^{-qt}
\le \max \left\{ \frac {\Psi}q, T(0)\right\}.
$$

\section{Conclusions}\label{sec:5} 

The analysis shows that the system cannot possibly disappear, as the origin is always
unstable. This is a good result in terms of ecological implications, and it is essentially
implicit in the model assumptions, namely the logistic growth of the prey, as the positive
eigenvalue of the Jacobian stems exactly from the susceptible prey reproduction equation.

The systems then have only two possible equilibria, related to each other via a transcritical
bifurcation, which occurs when the predators' mortality falls below the threshold $m^{\ddagger}$,
or its particular case $m^{\dagger}$ for the model of infected prey avoided by predators.
When it is above it, the systems settles at the prey-only equilibrium, with endemic disease.

The predator-free equilibrium
$E_1^c$ of the classical predator-prey case
can get destabilized by the disease presence, see the stability conditions for $E_1$
(\ref{E1_2_stab}) and (\ref{E1_stab}).
Note indeed that
the thresholds $m^{\dagger}$ and $m^{\ddagger}$ contain the epidemiological parameters $\beta$ and $\mu$,
while $m^*$ obviously does not, compare (\ref{m*}) with (\ref{E1_2_stab}) and (\ref{E1_stab}).
Note also that $m^{\dagger}=m^*$ for $\beta=0$, i.e. in the absence of the disease.
This destabilization never occurs in the classical case. This remark once more stresses the fact
that epidemics have also demographic consequences at the ecological level and therefore cannot be
easily neglected also in ecological investigations.

{\bf Acknowledgments:} This research has been partially supported by the project
``Metodi numerici in teoria delle popolazioni'' of the Dipartimento di Matematica ``Giuseppe Peano''.

\end{document}